\documentclass{article}[12pt]
\usepackage{amsmath,amsfonts,amssymb}
\usepackage[T2A]{fontenc}
\usepackage[cp866]{inputenc}
\usepackage[russian]{babel}
\newcommand{\la}{\lambda}
\newcommand{\en}{\enspace}

\begin{document}

\large

$$
$$

\bigskip

\begin{center}
\bf COMPLETE WAVE OPERATORS\\
IN NON-SELFADJOINT KATO MODEL OF\\
SMOOTH PERTURBATION THEORY
\end{center}

\bigskip

\begin{center}
\bf S.A.Stepin
\end{center}

\bigskip
\bigskip

{\bf Abstract.} Subject of the paper deals with the perturbation theory of 
 linear operators acting in Hilbert space. For a certain class of perturbations 
 the question is considered about existence of transformation operators 
 implementing linear similarity of perturbed and unperturbed operators. 
 In this context some results of complex analysis prove to be useful as well 
 as the relationship with the theory of operator semi-groups. 

\bigskip 

{\bf Keywords:} wave operators, scattering theory technique, Schroedinger operator

\bigskip

{\bf 2010 Mathematical Subject Classification:} 34L15, 34L25 

\bigskip

\begin{center}
\bf \S\,1. Introduction
\end{center}

\bigskip

The fundamental notion of relative smoothness was suggested by
T.Kato [1] within the framework of stationary scattering theory
technique elaborated therein for a certain class of not
necessarily selfadjoint operators. The goal of the approach worked
out in [1] was construction of wave operators intertwining
perturbed and unperturbed operators $\,T+V\,$ and $\,T.\,$ As a
matter of fact the original Kato paper deals with small
perturbations of the form $\,V=B^*A\,$ with operators $\,A\,$ and
$\,B\,$ smooth relatively to $\,T\,$ and $\,T^*\,$ respectively
(for precise definitions see Section 2). Besides it was assumed in
[1] that spectrum of $\,T\,$ is real and purely continuous while
operator $\,Q_0(\lambda)=A\,(T-\la I)^{-1}B^*\,$ has to be bounded
uniformly for non-real $\,\lambda.\,$ Wave operators constructed
in this context prove to be invertible and thus implement
similarity of initial and perturbed operators. Moreover the
statement of the problem also goes back to [1] about conditions
under which the approach in question admits a generalization to
the case of arbitrary (not necessarily small) perturbations of the
class under consideration. Note that the technique of factorizing
perturbation $\,V\,$ into the form $\,B^*A\,$ originated in [2].

Subsequent development of stationary scattering theory in the case
of relatively smooth perturbations basically dealt with
selfadjoint operators (see [3]). As regards applications of the
technique elaborated in [1] to non-selfadjoint operators $\,T+V\,$
when perturbation $\,V=B^*A\,$ is not assumed to be small, such
attempts essentially reduced (see e.g. [4],[5]) to construction of
local wave operators which supply with (and require) information
on the structure of perturbed operator $\,T+V\,$ in the subspace
corresponding to a certain segment of continuous spectrum. One of
the goals of of the present paper is to formulate conditions under
which stationary scattering theory approach in the framework of
non-selfadjoint Kato model context admits an extension to the case
of not necessarily small perturbations. Wave operators to be
constructed within this approach prove to possess completeness
property and thus realize the passage to spectral representation
of the perturbed operator implemented with regard for (and in
terms of) spectral continuity of the initial one.

Construction of local wave operators in the paper [5] was carried
out under the condition that operator $\,(I+Q_0(\la))^{-1}$ is
uniformly bounded in the neighborhood of the corresponding
continuous spectrum segment of operator $\,T.\,$ In the present
paper complete wave operators intertwining $\,T+V\,$ and $\,T,\,$
are going to be constructed provided that given resolvent
$\,R_V(\la)=(T+V-\la I)^{-1}\,$ the operator function
$\,Q_V(\lambda)=AR_V(\la)B^*\,$ is bounded uniformly for the
values of its argument separated away from discrete spectrum of
$\,T+V.\,$ As regards the question how to effectively verify this
condition it can be reduced to searching for a scalar analytic
function, most likely an appropriate Fredholm determinant, such
that all the singularities of operator function $\,Q_V(\la)\,$
after being multiplied by it become removable in $\,\mathbb
C_{\pm}.\,$ By virtue of resolvent identity
$$
R_V(\la)\,\,=\,\,R_0(\la)\,\,-\,\,R_0(\la)B^*(I+Q_0(\la))^{-1}AR_0(\la)
$$
condition suggested above is less restrictive compared to that
required in [5] while the class of perturbations $\,V=B^*A\,$
considered here is wider in due turn.

As an application we shall consider one-dimensional Schroedinger
operator $\,L\,=\,-d^2/dx^2+V(x)\,$ on half-axis $\,\mathbb R_+\,$
with complex potential $\,V(x)\,$ and Dirichlet boundary condition
at zero. Such an operator proves to be quite a simple and rather
capacious model which displays a number of effects typical for
perturbation theory in non-selfadjoint setting (see [6] and [7]).
Besides that Schroedinger operator with complex potential is known
(see [8]) to appear in the study of open quantum mechanical
systems with energy dissipation.

The approach elaborated in [1] reveals the following condition
\begin{equation*}
\int_0^{\infty}\!x\,|\,V(x)|\,dx\,\,<\,\,1
\end{equation*}
which guarantees similarity of $\,L\,$ and selfadjoint operator
$\,T=-d^2/dx^2\,$ corresponding to $\,V(x)\equiv 0.\,$ One of the
principal results of the present paper implies a criterion of the
similarity in question valid for the class of potentials
$\,V(x)\,$ possessing  a finite first momentum which extends and
supplements Kato sufficient condition [1].

\bigskip
\smallskip

\begin{center}
\bf \S\,2. Statement of basic results
\end{center}

\bigskip

As an unperturbed one consider closed operator $\,T\,$ with domain
$\,D(T)\,$ dense in Hilbert space $\,\cal H\,$ possessing purely
continuous spectrum $\,\sigma(T)\subset\mathbb R\,$ such that
resolvent $\,R_0(\la)=(T-\la I)^{-1}\,$ is analytic in $\,\mathbb
C_{\pm}.\,$ In accordance with [1] closed operator $\,A\,$ is said
to be smooth relative to $\,T\,$ if $\,AR_0(\la)\varphi\in{\rm
H}_2^{\pm}\,$ for arbitrary $\,\varphi\in\cal H,\,$ where $\,{\rm
H}_2^{\pm}\,$ denote Hardy classes corresponding to half-planes
$\,\mathbb C_{\pm},\,$ i.e.
$$
\sup_{\varepsilon>0}\,\int_{\mathbb R}\|AR_0(k\pm
i\varepsilon)\varphi\|^2\,dk\,<\,\infty\,.
$$
Moreover for almost all $\,k\in\mathbb R\,\,$ limits
$\,\,AR_0(k\pm i0)\varphi\,:=\,\lim\limits_{\varepsilon\downarrow
0}AR_0(k\pm i\varepsilon)\varphi\,\,$ exist (see [9]) and besides
$\,\, AR_0(k\pm i\varepsilon)\varphi\,\to\,AR_0(k\pm
i0)\varphi\,\,$ in mean square sense as $\,\varepsilon\downarrow
0.$

\smallskip

Along with $\,T\,$ introduce perturbed operator $\,L=T+V\,$ where
$\,V=B^*A\,$ so that $\,A\,$ is smooth relative to $\,T\,$ while
$\,B\,$ is smooth relative to $\,T^*.\,$ Below throughout the
following requirements are assumed to be fulfilled:

\bigskip

$(i)$ \, {\it resolvent $\,R_V(\la)\,$ is meromorphic in
$\,\mathbb C\setminus\sigma(T)\,$ so that its poles are just the
eigenvalues of $\,L\,$ with multiplicities taken into account;}

\medskip

$(ii)$ \, {\it discrete spectrum $\,\sigma_d(L)\,$ of operator
$\,L\,$ is a finite set while the complementary component
$\,\sigma_c(L)\,$ of its spectrum coincides with $\,\sigma(T);$}

\medskip

$(iii)$ \, {\it operator function $\,Q_V(\la)=AR_V(\la)B^*\,$ is
uniformly bounded in $\,\mathbb C_{\pm}\,$ provided that the
argument $\,\la\,$ is separated away from $\,\sigma_d(L).\,$}

\bigskip

For the sake of simplicity all the considerations will be carried
out here in the situation when $\,A\,$ and $\,B\,$ are bounded
operators and thus $\,D(L)=D(T).\,$ By means of contour integral
define Riesz projection
$$
P\,\,=\,-\,\frac1{2\pi i}\oint_{\Gamma}R_V(\la)\,d\la
$$
onto the linear span $\,{\cal H}_d:=P\cal H\,$ of root (eigen and
associated) vectors of operator $\,L\,$ parallel to the subspace
$\,{\cal H}_c:=(I-P)\cal H;\,$ here $\,\Gamma\,$ is appropriately
oriented closed contour separating mutually complementary spectral
components $\,\sigma_d(L)\,$ and $\,\sigma_c(L).\,$ Operator
$\,L\,$ can be decomposed with respect to direct sum
representation $\,{\cal H}={\cal H}_c\,\dot{+}\,{\cal H}_d\,$ in
the following (see [10]) sense
$$
PD(L)\subset D(L)\,,\quad L{\cal H}_c\subset{\cal H}_c\,,\quad
L{\cal H}_d\subset{\cal H}_d\,,
$$
while $\,\sigma(L|{\cal H}_c)=\sigma_c(L)\,$ and $\,\sigma(L|{\cal
H}_d)=\sigma_d(L).\,$ By $\,R_V^*(\la)\,$ and $\,{\cal
H}^*_d\,,\,{\cal H}^*_c\,$ let us denote the resolvent of operator
$\,L^*\,$ and its invariant subspaces corresponding to disjoint
spectral components $\,\sigma_d(L^*)\,$ and
$\,\sigma_c(L^*)\subset\mathbb R.\,$

\bigskip

{\bf Theorem 1.} {\it Assume that operator $\,T\,$ has purely
continuous spectrum $\,\sigma(T)\subset\mathbb R,\,$ operator
$\,V\,$ admits factorization $\,V=B^*A\,$ where $\,A\,$ and
$\,B\,$ are smooth relative to $\,T\,$ and $\,T^*\,$ respectively,
operator function $\,Q_0(\la)=AR_0(\la)B^*\,$ is bounded in
$\,\mathbb C_{\pm},\,$ while operator $\,L=T+V\,$ satisfies
conditions $\,(i)-(iii).\,$ Then there exist bounded stationary
wave operators $\,\,W\,$ and $\,Z\,$ being determined by their
bilinear forms {\,\rm:}
\begin{equation}
\big(W\varphi,\psi\big)\,\,=\,\,\big(\varphi,\psi\big)\,\,-\,\,\frac1{2\pi
i}\int_{\mathbb R}\big(AR_0(k+i0)\varphi,BR_V^*(k+i0)\psi\big)\,dk
\label{form1}
\end{equation}
for $\,\psi\in{\cal H}_c^*\,$ and $\,(W\varphi,\psi)=0\,$ if
$\,\psi\in{\cal H}_d^*\,;$
\begin{equation}
\big(Z\varphi,\psi\big)\,\,=\,\,\big(\varphi,\psi\big)\,\,+\,\,\frac1{2\pi
i}\int_{\mathbb
R}\big(AR_V(s+i0)\varphi,BR_0^*(s+i0)\psi\big)\,ds\,,
\label{form2}
\end{equation}
when $\,\varphi\in{\cal H}_c\,$ and $\,(Z\varphi,\psi)=0\,$ for
$\,\varphi\in{\cal H}_d.\,$ Besides $\,\,ZW=I\,\,$ and, moreover,
wave operators possess completeness property $\,\,W Z=I-P,\,\,$
where $\,P\,$ is a projection onto $\,{\cal H}_d=\ker Z\,$
parallel to $\,W{\cal H}={\cal H}_c.$  }

\bigskip

The proof of Theorem 1 is given in Sections 3 and 4 and it follows
the scheme of the approach elaborated in [1]. However in contrast
with [1] operators $\,A\,$ and $\,B\,$ are by no means smooth
relative to $\,L\,$ and $\,L^*\,$ respectively when
$\,\sigma_d(L)\ne\varnothing\,$ and such a crucial circumstance
proves to be of a considerable difficulty. Nevertheless the
construction of wave operators can be somewhat modified to treat
the setting in question due to the fact that
$\,AR_V(\la)\varphi\,$ and $\,BR_V^*(\la)\psi\in{\rm H}_2^{\pm}\,$
provided $\,\varphi\in{\cal H}_c\,$ and $\,\psi\in{\cal H}^*_c.\,$
To establish the fundamental properties of modified direct and
inverse wave operators we take advantage of the following
relationship of commutator type
$$ \frac1{2\pi
i}\,\lim_{\varepsilon\downarrow 0}\int_{\mathbb
R}R_0(k-i\varepsilon)\,[P,V]\,R_0(k+i\varepsilon)\,dk\,\,=\,\,P\,,
$$
where $\,[P,V]\,=\, PV-VP.\,$ Content of section 5 is summarized
by

\bigskip

{\bf Theorem 2.} {\it Under the assumptions of Theorem {\rm 1}
direct and inverse wave operators $\,W\,$ and $\,Z\,$ determined
by formulas {\rm (\ref{form1})} and {\rm (\ref{form2})} implement
similarity of $\,\,T\,\,$ and $\,\,L_c=L|{\cal H}_c\,$ so that
$\,\,L_c\,=\,WTZ.\,$ If additionally $\,T=T^*\,$ then one has the
following representations

\bigskip

$\quad\quad\quad\quad\quad\quad\quad\quad\quad
W\,\,=\,\,s$-$\!\lim\limits_{t\to\infty}\exp(itL)(I-P)\exp(-itT)\,,$

\medskip

$\quad\quad\quad\quad\quad\quad\quad\quad\quad
Z\,\,\,=\,\,\,s$-$\!\lim\limits_{t\to\infty}\exp(itT)\exp(-itL)(I-P)\,.$
}

\bigskip
\medskip

As an application of the approach developed here we shall consider
in Section 6 non-selfadjoint Schroedinger operator
$$
L\,\,=\,\,T\,+\,V\,\,=\,\,-\,d^2/dx^2\,+\,V(x)
$$
defined in the space $\,{\cal H}={\rm L}_2(\mathbb R_+)\,$ with
the domain $\,D(L)=\big\{y\in{\cal H}\!:y'$ absolutely continuous,
$y''\in{\cal H},$ $\,y(0)=0\big\},\,$ where potential $\,V(x)\,$
is a bounded complex-valued function. Provided that $\,V\in{\rm
L}_1(\mathbb R_+)\,$ equation
\begin{equation}
-\,y''\,\,+\,\,V(x)y\,\,=\,\,k^2y \label{form3}
\end{equation}
for any $\,k\in\mathbb C_+\,$ possesses Jost solution $\,e(x,k)\,$
which is asymptotically equivalent to $\,e^{ikx}\,$ at infinity.
In this situation spectrum of operator $\,L\,$ is known to consist
of continuous and discrete components
$$
\sigma_c(L)=\mathbb
R_+\,,\,\,\,\sigma_d(L)=\big\{\la=k^2\!:\,e(k)=0, k\in\mathbb
C_+\big\}\,,
$$
where $\,e(k):=e(0,k).\,$ The set of eigenvalues $\,\sigma_d(L)\,$
is bounded, at most countable and its accumulation points (if any)
belong to the half-axis $\,\mathbb R_+,\,$ whereas operator
$\,L\,$ has no eigenvalues embedded into continuous spectrum (see
[6]). Furthermore Jost function $\,e(k)\,$ is analytic in the open
half-plane $\,\mathbb C_+\,$ and admits extension by continuity to
$\,\mathbb R\setminus\{0\},\,$ while its real zeroes correspond to
distinguished points $\,\la=k^2\,$ of continuous spectrum
$\,\sigma_c(L)\,$ called spectral singularities (see [2]).

\bigskip

{\bf Theorem 3.} {\it Suppose that non-selfadjoint Schroedinger
operator
$$
L\,\,=\,\,T\,+\,V\,\,=\,\,-\,d^2/dx^2\,+\,V(x)
$$
defined in $\,{\cal H}={\rm L}_2(\mathbb R_+)\,$ by Dirichlet
boundary condition at zero and complex bounded potential
$\,V(x)\,$ satisfying condition
\begin{equation}
\int_0^{\infty}\!x\,|\,V(x)|\,dx\,\,<\,\,\infty\, \label{form4}
\end{equation}
has no spectral singularities. Then $\,L\,$ admits decomposition
with respect to direct sum $\,{\cal H}={\cal H}_c\,\dot{+}\,{\cal
H}_d\,$ generated by finite-dimensional Riesz projection $\,P\,$
onto the linear span $\,{\cal H}_d\,$ of eigen and associated
vectors of operator $\,L\,$ parallel to its continuous spectrum
invariant subspace $\,{\cal H}_c.\,$ Furthermore direct and
inverse stationary wave operators $\,W\,$ and $\,Z\,$ exist and
all the conclusions of Theorems {\,\rm 1} and {\rm 2} hold true.}

\bigskip

Absence of eigenvalues and spectral singularities in the context
of Theorem 3 leads to a criterion for similarity of operator
$\,L\,$ and selfadjoint operator $\,T\,$ and thus oblige $\,L\,$
to be spectral in the sense of Dunford. In the case when operator
$\,L\,$ does not have any spectral singularities results of papers
[6] and [7] enable one to produce transformation operators
intertwining $\,L\,$ and $\,T\,$ under rather restrictive
condition
\begin{equation*}
\int_0^{\infty}(1+x^2)\,|\,V(x)|\,dx\,\,<\,\,\infty\,.
\end{equation*}
However the corresponding procedure clarifies the relationship
between construction of wave operators and generalized
eigenfunction expansion problem for operator $\,L\,$ (cf. [11]).
Szoekefalvi-Nagy and Foias functional model technique was employed
by the author to construct wave operators given dissipative
Schroedinger operator with potential satisfying condition
(\ref{form4}) (see [12] and also references therein). The results
of the present paper in dissipative case were announced in [13]
(cf. [14]).

\bigskip
\smallskip

\begin{center}
\bf \S\,3. Completeness of wave operators
\end{center}

\bigskip

In what follows hypotheses of Theorem 1 are supposed to be
fulfilled.

\bigskip

{\bf Lemma 1.} {\it Given $\,\varphi\in{\cal H}_c\,$ and
$\,\psi\in{\cal H}^*_c\,$ vector functions $\,AR_V(\la)\varphi\,$
and $\,BR^*_V(\la)\psi\,$ belong to Hardy classes $\,{\rm
H}_2^{\pm}.$ }

\bigskip

In fact, provided that point $\,\la\not\in\sigma_d(L)\,$ is
contained inside $\,\Gamma\,$ by virtue of resolvent identity for
arbitrary $\,f\in{\cal H}\,$ one has
\begin{equation*}
R_V(\la)Pf\,\,=\,-\,\frac1{2\pi
i}\oint_{\Gamma}R_V(\la)R_V(\mu)f\,d\mu\,=
\,R_V(\la)f\,-\,\frac1{2\pi
i}\oint_{\Gamma}\frac{R_V(\mu)f}{\mu-\la}\,d\mu\,.
\end{equation*}
Hence given $\,\varphi=f-Pf\in{\cal H}_c\,$ vector function
$$
AR_V(\la)\varphi\,\,=\,\,\frac1{2\pi
i}\oint_{\Gamma}\frac{AR_V(\mu)f}{\mu-\la}\,d\mu
$$
admits analytic continuation to the points of discrete spectrum
$\,\sigma_d(L).\,$ Finally making use of relation
$\,AR_V(\la)\varphi=(I-Q_V(\la))AR_0(\la)\varphi,\,$ where
$\,AR_0(\la)\varphi\in{\rm H}_2^{\pm},\,$ and taking conditions
$\,(i)-(iii)\,$ into account we come to the required conclusion
$$
\sup_{\varepsilon>0}\,\int_{\mathbb R}\|AR_V(k\pm
i\varepsilon)\varphi\,\|^2\,dk\,<\,\infty\,.
$$
Similarly given $\,\psi\in{\cal H}^*_c\,$ vector function
$\,BR^*_V(\la)\psi\,$ extends analytically to the points
$\,\la\in\sigma_d(L^*)\,$ so that
$$
\sup_{\varepsilon>0}\,\int_{\mathbb R}\|BR^*_V(k\pm
i\varepsilon)\psi\,\|^2\,dk\,<\,\infty\,.
$$

\medskip

{\bf Corollary.} {\it Provided that $\,\varphi\in{\cal H}_c$ and
$\,\psi\in{\cal H}^*_c\,$ boundary values $\,AR_V(k\pm
i0)\varphi:=\lim\limits_{\varepsilon\downarrow 0}AR_V(k\pm
i\varepsilon)\varphi\,$ and $\,BR^*_V(k\pm i0)\psi:=$
$\lim\limits_{\varepsilon\downarrow 0}BR^*_V(k\pm
i\varepsilon)\psi\,$ are well-defined for almost all
$\,k\in\mathbb R\,$ and moreover $\,AR_V(k\pm
i\varepsilon)\varphi\,\to\,AR_V(k\pm i0)\varphi\,$ and
$\,BR^*_V(k\pm i\varepsilon)\psi\,\to\,BR^*_V(k\pm i0)\psi\,$ in
mean square sense as $\,\varepsilon\downarrow 0.$}

\bigskip

According to Lemma 1 and subsequent corollary the bilinear forms
(\ref{form1}) and (\ref{form2}) are bounded uniformly in
$\,\varphi,\psi\in\cal H,\,$ subject to normalization condition
$\,\|\varphi\|=\|\psi\|=1,\,$ and thus they determine properly
bounded operators $\,W\,$ and $\,Z.\,$ By corresponding
definitions one has $\,{\cal H}_d\subset\ker Z\,$ and $\,W{\cal
H}\perp{\cal H}^*_d\,$ where $\,{\cal H}^*_d\perp{\cal H}_c\,$ and
hence $\,W{\cal H}\subset{\cal H}_c.$

\bigskip

{\bf Statement 1.} {\it Wave operators $\,W\,$ and $\,Z\,$ possess
the completeness property{\,\rm :}
$$\,\,WZ\,=\,I-P.$$ }

\smallskip

{\bf Proof.} If $\,\varphi\in{\cal H}_d\subset\ker Z\,$ then
clearly $\,WZ\varphi=0.\,$ So one should verify that
$\,WZ\varphi=\varphi\,$ for arbitrary $\,\varphi\in{\cal H}_c.\,$
Since $\,W{\cal H}\subset{\cal H}_c\,$ one has
$\,\big(WZ\varphi,\psi\big)=0\,$ when $\,\psi\in{\cal H}^*_d.\,$
Hence it suffices to show that
$\,\big(WZ\varphi,\psi\big)=\big(\varphi,\psi\big)\,$ for
$\,\varphi\in{\cal H}_c\,$ and $\,\psi\in{\cal H}^*_c.\,$ To this
end making use of (\ref{form1}) and (\ref{form2}) consider
expression
\begin{multline*}
\big((W-I)(Z-I)\varphi,\psi\big)\,=\,\big((Z-I)\varphi,(W^*-I)\psi\big)\,\,=\\
=\,\,\frac1{2\pi i}\,\lim_{\varepsilon\downarrow 0}\,\int_{\mathbb
R}\big(AR_V(s+i\varepsilon)\varphi,BR_0^*(s+i\varepsilon)(W^*-I)\psi\big)\,ds\,,
\end{multline*}
where the integrand takes the form
\begin{multline*}
\big(AR_V(s+i\varepsilon)\varphi,BR_0^*(s+i\varepsilon)(W^*-I)\psi\big)\,=\,
\big((W-I)R_0(s-i\varepsilon)B^*AR_V(s+i\varepsilon)\varphi,\psi\big)\,=\\
=\,\,-\,\frac1{2\pi i}\,\int_{\mathbb
R}\big(AR_0(k+i0)R_0(s-i\varepsilon)B^*AR_V(s+i\varepsilon)\varphi,BR_V^*(k+i0)\psi\big)\,dk\,.
\end{multline*}
Here
$\,\,AR_0(k+i0)R_0(s-i\varepsilon)B^*\,=\,(k-s+i\varepsilon)^{-1}\big(Q_0(k+i0)-\,Q_0(s-i\varepsilon)\big)\,\,$
while boundary values $\,\,Q_0(k+i0)\,:=\,{\rm
s}\,$-$\lim\limits_{\varepsilon\downarrow 0}Q_0(k+
i\varepsilon)\,\,$ are well-defined almost everywhere and thus
\begin{multline*}
\big((W-I)(Z-I)\varphi,\psi\big)\,\,=\\
=\,\,\frac1{4\pi^2}\,\lim_{\varepsilon\downarrow
0}\,\bigg\{\int_{\mathbb R}dk\int_{\mathbb
R}(k-s+i\varepsilon)^{-1}\big(Q_0(k+i0)AR_V(s+i\varepsilon)\varphi,BR_V^*(k+i0)\psi\big)\,ds\,\,-\\
-\,\int_{\mathbb R}ds\int_{\mathbb
R}(k-s+i\varepsilon)^{-1}\big(Q_0(s-i\varepsilon)AR_V(s+i\varepsilon)\varphi,BR_V^*(k+i0)\psi\big)\,dk\bigg\}\,.
\end{multline*}
Note that $\,AR_V(\la)\varphi\,$ and $\,BR_V^*(\la)\psi\,$ belong
to $\,{\rm H}_2^{\pm}\,$ by virtue of Lemma 1 and one can apply
Cauchy integral formula to calculate transformations of Hilbert
type:
\begin{multline*}
\int_{\mathbb
R}(k-s+i\varepsilon)^{-1}\big(Q_0(k+i0)AR_V(s+i\varepsilon)\varphi,BR_V^*(k+i0)\psi\big)\,ds\,\,=\\
=\,-\,2\pi
i\,\big(AR_V(k+2i\varepsilon)\varphi,Q_0(k+i0)^*BR_V^*(k+i0)\psi\big)\,,
\end{multline*}
\begin{multline*}
\int_{\mathbb
R}(k-s+i\varepsilon)^{-1}\big(Q_0(s-i\varepsilon)AR_V(s+i\varepsilon)\varphi,BR_V^*(k+i0)\psi\big)\,dk\,\,=\\
=\,-\,2\pi
i\,\big(Q_0(s-i\varepsilon)AR_V(s+i\varepsilon)\varphi,BR_V^*(s+i\varepsilon)\psi\big)\,.
\end{multline*}
In order to pass to the limit in the right-hand side of equality
\begin{multline*}
\big((W-I)(Z-I)\varphi,\psi\big)\,\,=\\
=\,\,\frac1{2\pi i}\,\lim_{\varepsilon\downarrow
0}\,\bigg\{\int_{\mathbb R}\big(AR_V(k+2i\varepsilon)\varphi,Q_0(k+i0)^*BR_V^*(k+i0)\psi\big)\,dk\,\,-\\
-\,\int_{\mathbb
R}\big(Q_0(s-i\varepsilon)AR_V(s+i\varepsilon)\varphi,BR_V^*(s+i\varepsilon)\psi\big)\,ds\bigg\}\,,
\end{multline*}
we shall use corollary subsequent to Lemma 1 and besides the
existence of boundary values $\,Q_0(k\pm i0)\,$ almost everywhere
on the real axis. As a result (cf. [1]) it follows that
\begin{multline*}
\big((W-I)(Z-I)\varphi,\psi\big)\,\,=\\
=\,\,\frac1{2\pi i}\,\bigg\{\int_{\mathbb R}\big(Q_0(k+i0)AR_V(k+i0)\varphi,BR_V(k-i0)^*\psi\big)\,dk\,\,-\\
-\,\int_{\mathbb
R}\big(AR_V(s+i0)\varphi,Q_0(s-i0)^*BR_V(s-i0)^*\psi\big)\,ds\bigg\}\,\,=
\end{multline*}
\vspace{-0.6cm}
\begin{multline*}
=\,\,\frac1{2\pi i}\,\bigg\{\int_{\mathbb R}\big(AR_0(k+i0)\varphi,BR_V(k-i0)^*\psi\big)\,dk\,\,-\\
-\,\int_{\mathbb
R}\big(AR_V(s+i0)\varphi,BR_0(s-i0)^*\psi\big)\,ds\bigg\}\,,
\end{multline*}
since $\,\,Q_0(\la)AR_V(\la)\,=\,AR_0(\la)-AR_V(\la)\,\,$ and
$\,\,Q_0(\la)^*BR_V(\la)^*\,=\,BR_0(\la)^*-BR_V(\la)^*.\,$ Thus we
get
$$
\big((W-I)(Z-I)\varphi,\psi\big)\,\,=\,\,\big((I-W)\varphi,\psi\big)\,+\,\big((I-Z)\varphi,\psi\big)
$$
and hence $\,\big(WZ\varphi,\psi\big)=\big(\varphi,\psi\big)\,$
for arbitrary $\,\varphi\in{\cal H}_c,\,\psi\in{\cal H}^*_c.\,$

\bigskip

{\bf Corollary.} {\it Under the above assumptions one has
$\,W{\cal H}={\cal H}_c\,$ and $\,\,\ker Z={\cal H}_d.$ }

\bigskip
\smallskip

\begin{center}
\bf \S\,4. Inverse wave operator
\end{center}

\bigskip

Along with projection $\,P\,$ introduce the complementary
projection $\,\widetilde{P}=I-P\,$ onto the subspace $\,{\cal
H}_c\,$ parallel to the linear span $\,{\cal H}_d\,$ of eigen and
associated vectors of operator $\,L.\,$

\bigskip

{\bf Lemma 2.} {\it Under the hypotheses $\,(i)-(iii)\,$ operator
function $\,\widetilde{Q}(\la)=AR_V(\la)\widetilde{P}B^*\,$ is
analytic and bounded in $\,\mathbb C_{\pm}.$}

\bigskip

Indeed, due to conditions $\,(i)-(ii)\,$ operator function
$\,AR_V(\la)\widetilde{P}B^*\,$ being extended by the expression
$$
\frac1{2\pi
i}\oint_{\Gamma}\,\frac{AR_V(\zeta)\widetilde{P}B^*}{\zeta-\la}\,\,d\zeta
$$
to singular points $\,\la\in\sigma_d(L)\,$ becomes analytic in
$\,\mathbb C_{\pm}\,$ and bounded in certain neighborhood of
$\,\sigma_d(L).\,$ Provided that $\,\la\in\mathbb C_{\pm}\,$ is
separated away from $\,\sigma_d(L)\,$ operator function
$$
\widetilde{Q}(\la)\,+\,AR_V(\la)PB^*\,=\,AR_V(\la)B^*
$$
is uniformly bounded by condition $\,(iii).\,$ Hence it suffices
to show that operator function $\,AR_V(\la)PB^*\,$ is bounded
uniformly when its argument $\,\la\in\mathbb C_{\pm}\,$ does not
approach $\,\sigma_d(L).\,$ In this context according to [15]
operator $\,P\,$ proves to be the sum of Riesz projections onto
root subspaces corresponding to the points of discrete spectrum
$\,\sigma_d(L).\,$ Given an eigenvalue $\,\mu\,$ of operator
$\,L\,$ with geometric multiplicity $\,n\,$ the associated Riesz
projection is of the form
$$
P_{\mu}\,\,=\,\,\sum_{j=1}^n\,\big\{\,(\,\cdot\,,g_j^{(m_j-1)}\big)f_j^{(0)}\,+\,\ldots\,+
\,\,(\,\cdot\,,g_j^{(0)}\big)f_j^{(m_j-1)}\big\}\,.
$$

\vspace{-0.3cm}

\noindent Eigenvectors $\,f_j^{(0)},\,j=1,\ldots,n,\,$ span the
eigensubspace $\,\ker(L-\mu I)\,$ of dimension $\,n,\,$ while
$\,f_j^{(1)},\ldots, f_j^{(m_j-1)}\,$ is a Jordan chain of
associated vectors adjoint to $\,f_j^{(0)},\,$ so that
$\,m=m_1+\ldots+m_n\,$ is the algebraic multiplicity of $\,\mu.\,$
In due turn $\,g_j^{(0)},g_j^{(1)},\ldots, g_j^{(m_j-1)}\,$ is the
chain of eigen and associated vectors of operator $\,L^*\,$
corresponding to its eigenvalue $\,\overline\mu\,$ subject to
normalization condition $\,(f_j^{(p)},g_j^{(q)})\,=\,1\,$ where
$\,p+q=m_j-1.\,$ Since $\,P_{\mu}\,$ reduces to the sum of rank 1
operators
$\,\,\widehat{P}_{\mu}=\big(\,\cdot\,,g^{(q)}_j\big)f_j^{(p)}\,$
it suffices to consider operator function
$\,AR_V(\la)\widehat{P}_{\mu}B^*\,$ which by virtue of the
estimate
$$
\|AR_V(\la)\widehat{P}_{\mu}B^*\|\,\,\leqslant\,\,\sum_{r=0}^p\,|\la-\mu\,|^{r-p-1}\|Af_j^{(r)}\|\,\|Bg_j^{(q)}\|\,
$$
proves to be uniformly bounded in $\,\mathbb C_{\pm}\,$ outside
some neighborhood of $\,\sigma_d(L).\,$

\bigskip

{\bf Corollary.} \, {\it Boundary values $\,\,\,\widetilde{Q}(k\pm
i0)\,\,:=\,\,{\rm s}\,$-$\lim\limits_{\varepsilon\downarrow
0}\widetilde{Q}(k\pm i\varepsilon)\,\,\,$ exist for almost all
$\,k\in\mathbb R.\,$ }

\bigskip

{\bf Statement 2.} {\it Operator $\,Z\,$ is the left inverse to
the direct wave operator $\,W.$ }

\bigskip

{\bf Proof.} For arbitrary $\,\varphi,\psi\in\cal H\,$ one has
\begin{multline*}
\big(ZW\varphi,\psi\big)\,=\,\big(W\varphi,Z^*\psi\big)\,\,=\\
=\,\,\big(\varphi,Z^*\psi\big)\,\,-\,\,\frac1{2\pi
i}\,\lim_{\varepsilon\downarrow 0}\,\int_{\mathbb
R}\big(AR_0(s+i\varepsilon)\varphi,BR^*_V(s+i\varepsilon)Z^*\psi\big)\,ds\,,
\end{multline*}
where $\,Z^*\psi\in{\cal H}^*_c,\,$ so that the integrand in
virtue of (\ref{form2}) takes the form
\begin{multline*}
\big(ZR_V(s-i\varepsilon)VR_0(s+i\varepsilon)\varphi,\psi\big)\,\,=\,\,
\big(\widetilde{P}R_V(s-i\varepsilon)VR_0(s+i\varepsilon)\varphi,\psi\big)\,\,+\\
+\,\,\frac1{2\pi i}\int_{\mathbb
R}\big(AR_V(k+i0)\widetilde{P}R_V(s-i\varepsilon)VR_0(s+i\varepsilon)\varphi,BR_0(k-i0)^*\psi\big)\,dk\,.
\end{multline*}
Thus we get
\begin{multline*}
\big(ZW\varphi,\psi\big)\,=\,\big(Z\varphi,\psi\big)\,-\,
\frac1{2\pi i}\,\lim_{\varepsilon\downarrow 0}\,\int_{\mathbb R}
\big(\widetilde{P}R_V(s-i\varepsilon)VR_0(s+i\varepsilon)\varphi,\psi\big)\,ds\,\,\,+\\
+\,\,\frac1{4\pi^2}\,\lim_{\varepsilon\downarrow 0}\,\int_{\mathbb
R}ds\int_{\mathbb
R}\big(AR_V(k+i0)\widetilde{P}R_V(s-i\varepsilon)B^*AR_0(s+i\varepsilon)\varphi,BR_0(k-i0)^*\psi\big)\,dk\,
\end{multline*}
and according to resolvent identity
$$
AR_V(k+i0)\widetilde{P}R_V(s-i\varepsilon)B^*\,\,=\,\,
(k-s+i\varepsilon)^{-1}\big(\widetilde{Q}(k+i0)\,-\,\widetilde{Q}(s-i\varepsilon)\big)\,
$$
let us split the above iterated integral into two summands
\begin{multline*}
\int_{\mathbb R}ds\int_{\mathbb
R}\big(AR_V(k+i0)\widetilde{P}R_V(s-i\varepsilon)B^*AR_0(s+i\varepsilon)\varphi,BR_0(k-i0)^*\psi\big)\,dk\,\,=\\
=\,\,\int_{\mathbb R}dk\int_{\mathbb
R}(k-s+i\varepsilon)^{-1}\big(AR_0(s+i\varepsilon)\varphi,\widetilde{Q}(k+i0)^*BR_0(k-i0)^*\psi\big)\,ds\,\,\,-\\
-\,\,\int_{\mathbb R}ds\int_{\mathbb
R}(k-s+i\varepsilon)^{-1}\big(\widetilde{Q}(s-i\varepsilon)AR_0(s+i\varepsilon)\varphi,BR_0(k-i0)^*\psi\big)\,dk\,.
\end{multline*}

Making use of Cauchy integral formula (cf. the proof of Statement
1) one can compute Hilbert-type transforms:
\begin{multline*}
\int_{\mathbb
R}(k-s+i\varepsilon)^{-1}\big(AR_0(s+i\varepsilon)\varphi,\widetilde{Q}(k+i0)^*BR_0(k-i0)^*\psi\big)\,ds\,\,=\\
=\,\,-\,2\pi
i\,\big(AR_0(k+2i\varepsilon)\varphi,\widetilde{Q}(k+i0)^*BR_0(k-i0)^*\psi\big)\,,
\end{multline*}
\vspace{-0.7cm}
\begin{multline*}
\int_{\mathbb
R}(k-s+i\varepsilon)^{-1}\big(\widetilde{Q}(s-i\varepsilon)AR_0(s+i\varepsilon)\varphi,BR_0(k-i0)^*\psi\big)\,dk\,\,=\\
=\,\,-\,2\pi
i\,\big(\widetilde{Q}(s-i\varepsilon)AR_0(s+i\varepsilon)\varphi,BR_0(s-i\varepsilon)^*\psi\big)\,;
\end{multline*}
now passing in the expression for $\,\big(ZW\varphi,\psi\big)\,$
to the limit as $\,\varepsilon\downarrow 0\,$ leads to the
following representation
\begin{multline*}
\big(ZW\varphi,\psi\big)\,=\,\big(Z\varphi,\psi\big)\,-\,
\frac1{2\pi i}\,\int_{\mathbb R}
\big(AR_0(s+i0)\varphi,BR_V(s-i0)^*\widetilde{P}^*\psi\big)\,ds\,\,\,+\\
+\,\,\frac1{2\pi i}\,\int_{\mathbb R}
\big(\widetilde{Q}(k+i0)AR_0(k+i0)\varphi,BR_0(k-i0)^*\psi\big)\,dk\,\,\,-\\
-\,\,\frac1{2\pi i}\,\int_{\mathbb R}
\big(AR_0(s+i0)\varphi,\widetilde{Q}(s-i0)^*BR_0(s-i0)^*\psi\big)\,ds\,.
\end{multline*}
Due to equalities
$\,\,\widetilde{Q}(\la)AR_0(\la)=A\widetilde{P}R_0(\la)\,-\,AR_V(\la)\widetilde{P}\,$
and
$\,\widetilde{Q}(\la)^*BR_0(\la)^*=\,B\widetilde{P}^*R_0(\la)^*-BR_V(\la)^*\widetilde{P}^*\,$
for arbitrary $\,\varphi,\psi\in{\cal H}\,$ and almost all
$\,k\in\mathbb R\,$ in virtue of Lemmas 1 and 2 there exist
boundary values
$\,\,A\widetilde{P}R_0(k+i0)\varphi,\,B\widetilde{P}^*R_0(k-i0)^*\psi\in{\rm
L}_2(\mathbb R).\,$ This fact enables one to reduce bilinear form
of operator $\,ZW\,$ to the following expression
\begin{gather*}
\big(ZW\varphi,\psi\big)\,\,=\,\,\big((I-P)\varphi,\psi\big)\,\,+\,\,
\frac1{2\pi i}\,\int_{\mathbb R}
\big(A\widetilde{P}R_0(k+i0)\varphi,BR_0(k-i0)^*\psi\big)\,dk\,\,-\\
-\,\,\,\frac1{2\pi i}\,\int_{\mathbb R}
\big(AR_0(k+i0)\varphi,B\widetilde{P}^*R_0(k-i0)^*\psi\big)\,dk\,\,=\,\,\big((I-P)\varphi,\psi\big)\,\,+\\
+\,\,\frac1{2\pi i}\,\lim_{\varepsilon\downarrow 0}\int_{\mathbb
R}\big[\big(R_0(k-i\varepsilon)PVR_0(k+i\varepsilon)\varphi,\psi\big)-
\big(R_0(k-i\varepsilon)VPR_0(k+i\varepsilon)\varphi,\psi\big)\big]\,dk\,,
\end{gather*}
where
\begin{multline*}
R_0(k-i\varepsilon)PVR_0(k+i\varepsilon)\,-\,R_0(k-i\varepsilon)PVR_0(k+i\varepsilon)\,\,=\\
=\,\,PR_0(k+i\varepsilon)\,-\,R_0(k-i\varepsilon)P\,-\,2i\varepsilon
R_0(k-i\varepsilon)PR_0(k+i\varepsilon)\,.
\end{multline*}

In order to complete the proof we only need to establish the
relationship
\begin{multline}
\frac1{2\pi i}\,\lim_{\varepsilon\downarrow 0}\int_{\mathbb
R}\big[\big(PR_0(k+i\varepsilon)\varphi,\psi\big)\,-\,\big(R_0(k-i\varepsilon)P\varphi,\psi\big)\,\,-\\
-\,\,2i\varepsilon
\big(R_0(k-i\varepsilon)PR_0(k+i\varepsilon)\varphi,\psi\big)\big]\,dk\,\,=\,\,\big(P\varphi,\psi\big)\,.
\label{form6}
\end{multline}
Since operator $\,P\,$ turns out to be the sum of Riesz
projections $\,P_{\mu}\,$ corresponding to points of discrete
spectrum $\,\mu\in\sigma_d(L)\,$ it suffices to verify relation
(\ref{form6}) separately for each elementary rank 1 summand
$\,\widehat{P}_{\mu}=\big(\,\cdot\,,g\big)f,\,$ where
$\,f=f_j^{(p)},\,g=g_j^{(q)},$ $\,p+q=m_j-1\,$ and
$\,\big(f,g\big)=1.\,$ Note that
$$
\big(\widehat{P}_{\mu}R_0(k+i\varepsilon)\varphi,\psi\big)\,-\,\big(R_0(k-i\varepsilon)\widehat{P}_{\mu}\varphi,\psi\big)\,-
\,2i\varepsilon
\big(R_0(k-i\varepsilon)\widehat{P}_{\mu}R_0(k+i\varepsilon)\varphi,\psi\big)\,\,=\\
$$
\vspace{-0.9cm}
\begin{multline*}
=\,\,\big(R_0(k+i\varepsilon)\varphi,g\big)\big(f,\psi\big)\,-\,\big(\varphi,g\big)\big(R_0(k-i\varepsilon)f,\psi\big)\,\,-\\
-\,\,2i\varepsilon\big(R_0(k+i\varepsilon)\varphi,g\big)\big(R_0(k-i\varepsilon)f,\psi\big)\,.
\end{multline*}
Taking representations $\,\,\displaystyle{f\,=\,-\sum_{r=0}^p
R_0(\mu)^{p-r+1}Vf_j^{(r)},\,\,g\,=\,-\sum_{s=0}^qR^*_0(\overline{\mu})^{q-s+1}V^*g_j^{(s)}}\,\,$
into account according to resolvent identity we get
\begin{eqnarray*}
\big(R_0(k+i\varepsilon)\varphi,g\big)&=&-\,\,\frac{\big(\varphi,g\big)}{k-\mu+i\varepsilon}\,\,-\,\,
\frac{\xi(k+i\varepsilon)}{k-\mu+i\varepsilon}\,,\\
\big(R_0(k-i\varepsilon)f,\psi\big)&=&-\,\,\frac{\big(f,\psi\big)}{k-\mu-i\varepsilon}\,\,-\,\,
\frac{\eta(k-i\varepsilon)}{k-\mu-i\varepsilon}\,,
\end{eqnarray*}
where
$\,\,\displaystyle{\xi(z)\,=\sum_{s=0}^q\big(R_0(\mu)^{q-s}\varphi,R_0(z)^*V^*g_j^{(s)}\big),\,\,
\eta(z)\,=\sum_{r=0}^p\big(R_0(z)Vf_j^{(r)},R_0^*(\overline{\mu})^{p-r}\psi\big)}\,\,$
belong to Hardy classes $\,{\rm H}_2^{\pm}\,$ so that
\begin{multline*}
\big(\widehat{P}_{\mu}R_0(k+i\varepsilon)\varphi,\psi\big)\,-\,
\big(R_0(k-i\varepsilon)\widehat{P}_{\mu}\varphi,\psi\big)\,-\,2i\varepsilon
\big(R_0(k-i\varepsilon)\widehat{P}_{\mu}R_0(k+i\varepsilon)\varphi,\psi\big)\,\,=\\
=\,\,\frac{\big(\varphi,g\big)}{k-\mu+i\varepsilon}\,\eta(k-i\varepsilon)\,\,-\,\,\frac{\big(f,\psi\big)}{k-\mu-i\varepsilon}\,
\xi(k+i\varepsilon)\,\,
-\,\,\frac{2i\varepsilon}{(k-\mu)^2+\varepsilon^2}\,\xi(k+i\varepsilon)\,\eta(k-i\varepsilon)\,.
\end{multline*}
In the case when $\,{\rm Im}\,\mu>0\,$ by Cauchy integral formula
one has
$$
\frac1{2\pi i}\,\int_{\mathbb
R}\frac1{k-\mu-i\varepsilon}\,\xi(k+i\varepsilon)\,dk\,\,=\,\,
\xi(\mu+2i\varepsilon)\,\,\to\,\,-\,\big(\varphi,g\big)\,,\enspace
\varepsilon\downarrow 0\,,
$$
while for sufficiently small $\,\varepsilon>0\,$
$$
\int_{\mathbb
R}\frac1{k-\mu+i\varepsilon}\,\eta(k-i\varepsilon)\,dk\,\,=\,\,0\,.
$$
Besides the following estimate is valid
\begin{multline*}
\left|\int_{\mathbb
R}\frac1{(k-\mu)^2+\varepsilon^2}\,\xi(k+i\varepsilon)\,\eta(k-i\varepsilon)\,dk\,\right|\,\,\,\leqslant\\
\leqslant\,\,\,\frac1{(\,{\rm
Im}\,\mu)^2-\varepsilon^2}\,\bigg(\int_{\mathbb
R}|\,\xi(k+i\varepsilon)|^2dk\bigg)^{1/2}\bigg(\int_{\mathbb
R}|\,\eta(k-i\varepsilon)|^2dk\bigg)^{1/2},
\end{multline*}
where the right-hand side is bounded uniformly provided
$\,\varepsilon>0\,$ is small enough. Thus given eigenvalue
$\,\mu\in\mathbb C_+\,$ the limiting relationship (\ref{form6}) is
established for $\,\widehat{P}_{\mu}\,$ and similarly it can be
derived if $\,\mu\in\sigma_d(L)\cap\mathbb C_-.\,$ As regards the
case of real $\,\mu\in\sigma_d(L)\,$ one should take into account
analyticity of resolvent $\,R_0(\la)\,$ at such points to carry
out passage to the limit as $\,\varepsilon\downarrow 0\,$
\begin{eqnarray*}
\frac1{2\pi i}\,\int_{\mathbb
R}\frac1{k-\mu-i\varepsilon}\,\xi(k+i\varepsilon)\,dk\,\,&=&\,\,
\xi(\mu+2i\varepsilon)\,\,\to\,\,-\,\big(\varphi,g\big)\,,\\
\frac1{2\pi i}\,\int_{\mathbb
R}\frac1{k-\mu+i\varepsilon}\,\eta(k-i\varepsilon)\,dk\,\,&=&\,\,
-\,\eta(\mu-2i\varepsilon)\,\,\to\,\,\big(f,\psi\big)\,.
\end{eqnarray*}
By the same argument we get
$$
\lim_{\varepsilon\downarrow
0}\,\frac{\varepsilon}{\pi}\int_{\mathbb
R}\frac1{(k-\mu)^2+\varepsilon^2}\,\xi(k+i\varepsilon)\,\eta(k-i\varepsilon)\,dk\,\,=\,\,
\big(\varphi,g\big)\,\big(f,\psi\big)\,,
$$
and so the above calculations imply relationship (\ref{form6}) for
rank 1 operator $\widehat{P}_{\mu}$ in question.

\bigskip
\medskip

\begin{center}
\bf \S\,5. Non-stationary wave operators
\end{center}

\bigskip

Let us first verify that wave operators constructed in Theorem 1
and possessing completeness property implement similarity of
operator $\,T\,$ and the restriction of operator $\,L\,$ to the
subspace $\,{\cal H}_c.$

\bigskip

{\bf Statement 3.} {\it Wave operator $\,W\,$ intertwines the
resolvents of perturbed and unperturbed operators $\,L\,$ and
$\,T$\,{\rm :} $\,\, R_V(\la)\,W\,=\,WR_0(\la).$ }

\bigskip

{\bf Proof.} To be specific we shall consider the case $\,{\rm
Im}\,\la>0.\,$ Given $\,\psi\in{\cal H}^*_d\,$ for arbitrary
$\,\varphi\in\cal H\,$ one has
$\,\big(WR_0(\la)\varphi,\psi\big)=0\,\,$ and
$\,\big(R_V(\la)W\varphi,\psi\big)$
$=\big(W\varphi,R_V(\la)^*\psi\big)$ $=0\,$ since
$\,R_V(\la)^*\psi\in{\cal H}^*_d.\,$

If $\,\psi\in{\cal H}^*_c\,$ then $\,BR_V^*(\la)\psi\in{\rm
H}_2^{\pm}\,$ by virtue of Lemma 1. According to (\ref{form1}) the
following representation
$$
\big((W-I)R_0(\la)\varphi,\psi\big)\,\,=\,\,-\,\frac1{2\pi
i}\int_{\mathbb
R}\big(AR_0(k+i0)R_0(\la)\varphi,BR_V^*(k+i0)\psi\big)\,dk\,
$$
is valid and due to resolvent identity
$$
AR_0(k+i0)R_0(\la)\varphi\,\,=\,\,(k-\la)^{-1}\big(AR_0(k+i0)\varphi-AR_0(\la)\varphi\big)
$$
we hence get
$$
\big((W-I)R_0(\la)\varphi,\psi\big)\,\,=\,\,-\,\frac1{2\pi
i}\int_{\mathbb
R}(k-\la)^{-1}\big(AR_0(k+i0)\varphi,BR_V^*(k+i0)\psi\big)\,dk\,
$$
because
$$
\int_{\mathbb
R}(k-\la)^{-1}\big(AR_0(\la)\varphi,BR_V^*(k+i0)\psi\big)\,dk\,\,=\,\,0\,.
$$

Similarly noting that $\,R_V^*(\la)\psi\in{\cal H}_c^*\,$ one can
carry out the following transformation
\begin{multline*}
\big(R_V(\la)(W-I)\varphi,\psi\big)\,\,=\,\,\big((W-I)\varphi,R_V(\la)^*\psi\big)\,\,=\\
=\,\,\frac1{2\pi i}\int_{\mathbb
R}(k-\la)^{-1}\big(AR_0(k+i0)\varphi,BR_V(\la)^*\psi\big)\,dk\,\,\,-\\
-\,\,\frac1{2\pi i}\int_{\mathbb
R}(k-\la)^{-1}\big(AR_0(k+i0)\varphi,BR_V^*(k+i0)\psi\big)\,dk\,,
\end{multline*}
where the first summand on the right-hand side by virtue of Cauchy
integral formula reduces to the expression
$$
\big(AR_0(\la)\varphi,BR_V(\la)^*\psi\big)\,=\,\big(R_V(\la)VR_0(\la)\varphi,\psi\big)\,=\,
\big(R_0(\la)\varphi,\psi\big)-\big(R_V(\la)\varphi,\psi\big)\,.
$$
Finally we get the equality
$$
\big(R_V(\la)(W-I)\varphi,\psi\big)\,=\,\big(R_0(\la)\varphi,\psi\big)-\big(R_V(\la)\varphi,\psi\big)+
\big((W-I)R_0(\la)\varphi,\psi\big)\,
$$
and thus for arbitrary vectors $\,\varphi,\psi\in\cal H\,$ the
required relationship
$\,\big(R_V(\la)W\varphi,\psi\big)=\big(WR_0(\la)\varphi,\psi\big)\,$
holds.

\bigskip

{\bf Corollary.} {\it Under hypotheses of Theorem {\rm 1} one has
$\,\,L_c=\,WTZ.$}

\bigskip

In the case of selfadjoint unperturbed operator $\,T\,$ along with
unitary group $\,\,U_0(t):=\exp\big(itT\big)\,$ we put into
context a one-parameter operator family

\medskip

$\quad\quad\quad\quad\quad\quad\quad\quad\quad
\exp\big(itL\big)\,\,=\,\,s$-$\!\lim\limits_{n\to\infty}\big(I-itL/n\big)^{-n}.$

\bigskip

{\bf Lemma 3.} {\it For arbitrary $\,t\in\mathbb R\,$ exponents
$\,\,\exp\big(itT\big)\,$ and $\,\,\exp\big(itL\big)\,$ satisfy
the equation
$$
U_V(t)\,:=\,\exp\big(itL\big)\widetilde{P}\,\,=\,\,W\,U_0(t)\,Z.
$$
}

\vspace{-0.2cm}

Indeed, the intertwining relation for the resolvents
$\,R_V(\la)\,$ and $\,R_0(\la)\,$ (Statement 3) implies the
equality
$$
\big(I-itL/n\big)^{-n}W\,\,=\,\,W\,\big(I-itT/n\big)^{-n},
$$
which by the use of completeness property for wave operators
(Statement 1) can be rewritten in the form
$$
\big(I-itL/n\big)^{-n}(I-P)\,\,=\,\,W\,\big(I-itT/n\big)^{-n}Z.
$$
It thus follows that for arbitrary $\,\varphi\in{\cal H}\,$ the
limit
$$
\lim\limits_{n\to\infty}\big(I-itL/n\big)^{-n}\widetilde{P}\,\varphi\,\,=\,\,WU_0(t)Z\,\varphi
$$
exists, i.e. $\,U_V(t)=WU_0(t)Z.\,$ Remark that exponent
$\,\exp\big(itL\big)\,$ is well-defined on the whole space $\,\cal
H:\,$ namely if $\,\varphi\in{\cal H}_c\,$ then
$\,\,\exp\big(itL\big)\varphi\,=\,WU_0(t)Z\varphi,\,$ while one
calculates the value $\,\exp\big(itL\big)\varphi\,$ for
$\,\varphi\in{\cal H}_d\,$ taking notice of the action formulas
$$
\exp\big(itL\big)\,f_j^{(r)}\,\,=\,\,\,e^{it\mu}\,\sum_{s=0}^r\,\frac{(it)^{r-s}}{(r-s)!}\,\,f_j^{(s)}
$$
in the Jordan basis
$\,\{f_j^{(0)},f_j^{(1)},\ldots,f_j^{(m_j-1)}\}_{j=1}^n\,$ of the
root subspace corresponding to eigenvalue $\,\mu\,$ of operator
$\,L.\,$

\bigskip

{\bf Corollary.} {\it One-parameter operator family $\,U_V(t)\,$
possesses semigroup property $\,\,U_V(t+s)\,=\,U_V(t)\,U_V(s)\,\,$
and is related to operator function $\,R_V(\la)\widetilde{P}\,$ by
means of Laplace transform. }

\bigskip

{\bf Statement 4.} {\it If additionally to the hypotheses of
Theorem {\rm 1} unperturbed operator $\,T\,$ is assumed to be
selfadjoint then direct and inverse wave operators $\,W\,$ and
$\,Z\,$ admit representations

\medskip

$\quad\quad\quad\quad\quad
W\,=\,s$-$\!\lim\limits_{t\to\infty}U_V(t)\,U_0(-t)\,, \quad
Z\,=\,s$-$\!\lim\limits_{t\to\infty}U_0(t)\,U_V(-t)\,.$ }

\bigskip
\smallskip

{\bf Proof\,\,} follows the scheme of the approach suggested in
[1]. Making use of Lemmas 1 and 3 and taking into account
smoothness of operator $\,A\,$ relative to $\,T\,$ we come to a
conclusion that for $\,\varphi\in{\cal H}\,$ and $\,\psi\in{\cal
H}^*_c\,$ images of vector functions $\,\,AR_0(k+i0)\varphi\,\,$
and $\,\,BR^*_V(k+i0)\psi\,\,$ under the action of
Fourier-Plancherel transform are given by the expressions
$\,\,i\,\sqrt{2\pi}\,\theta(t)AU_0(-t)\varphi\,\,$ and
$\,\,i\,\sqrt{2\pi}\,\theta(t)BU^*_V(-t)\psi,\,\,$ where
$\,\theta(t)\,$ is Heaviside step function and
$\,\,U_V^*(-t)\,=\,\exp\big(-itL^*\big)\widetilde{P}^*=\,U_V(t)^*.\,\,$
Therefore by virtue of Parseval identity one has
$$
\int_{\mathbb
R}\big(AR_0(k+i0)\varphi,BR_V^*(k+i0)\psi\big)\,dk\,\,=\,\,
2\pi\int_0^{\infty}\!\!\big(AU_0(-t)\varphi,BU^*_V(-t)\psi\big)\,dt\,.
$$
Due to this fact and with regard for intertwining relationship
$\,\,U_V(s)W=WU_0(s)\,\,$ and semigroup property
$\,\,U_V(t+s)=U_V(t)U_V(s)\,\,$ (corollary subsequent to Lemma 3)
bilinear form of wave operator $\,W\,$ is reduced to the
expression
\begin{multline*}
\big(W\varphi,\psi\big)\,\,=\,\,\big(WU_0(-s)\varphi,U_V^*(-s)\psi\big)\,\,\,=\\
=\,\,\,\big(U_V(s)U_0(-s)\varphi,\psi\big)\,\,+\,\,\,i\int_s^{\infty}\!\!\big(AU_0(-t)\varphi,BU^*_V(-t)\psi\big)\,dt\,
\end{multline*}
which implies the estimate
$$
\big|\big(W\varphi-U_V(s)U_0(-s)\varphi,\psi\big)\big|\,\leqslant\,
\bigg(\!\int_s^{\infty}\!\!\|AU_0(-t)\varphi\|^2\,dt\bigg)^{1/2}\!\!
\bigg(\!\int_0^{\infty}\!\!\|BU^*_V(-t)\psi\|^2\,dt\bigg)^{1/2}\!,
$$
where
$$
\int_0^{\infty}\!\!\|BU^*_V(-t)\psi\|^2\,dt\,\,\,=\,\,\,\frac1{2\pi}\int_{\mathbb
R}\|BR^*_V(k+i0)\psi\|^2\,dk\,.
$$
Given $\,\psi\in{\cal H}^*_c\,$ by virtue of Lemma 1 and the
subsequent corollary a constant $\,C>0\,$ exists such that the
inequality
$$
\big|\,\big(W\varphi-U_V(s)U_0(-s)\varphi,\psi\big)\,\big|\,\,\leqslant\,\,C\,
\bigg(\int_s^{\infty}\!\!\|AU_0(-t)\varphi\|^2\,dt\bigg)^{1/2}\!\!\|\psi\|
$$
is valid for arbitrary $\,\varphi\in{\cal H}.\,$ If however
$\,\psi\in{\cal H}^*_d\,$ then $\,\big(W\varphi,\psi\big)=0\,$ by
the definition of $\,W\,$ and moreover
$\,\big(U_V(s)U_0(-s)\varphi,\psi\big)=0\,$ because
$\,U_V(s)U_0(-s)\varphi\in{\cal H}_c\perp{\cal H}^*_d.\,$ Thus one
has
$$
\big\|\,W\varphi-U_V(s)U_0(-s)\varphi\,\big\|\,\,\leqslant\,\,C\,
\bigg(\int_s^{\infty}\!\!\|AU_0(-t)\varphi\|^2\,dt\bigg)^{1/2},
$$
where $\,AU_0(-t)\varphi\in{\rm L}_2(\mathbb R_+)\,$ for arbitrary
$\,\varphi\in{\cal H}\,$ and therefore wave operator $\,W\,$
coincides with the strong limit of $\,U_V(s)\,U_0(-s)\,$ as
$\,s\to\infty.\,$

Similarly with Lemmas 1 and 3 taken into account and in view of
smoothness of operator $\,B\,$ relative to $\,T^*\,$ the equality
$$
\int_{\mathbb
R}\big(AR_V(k+i0)\varphi,BR_0(k+i0)\psi\big)\,dk\,\,=\,\,
2\pi\int_0^{\infty}\!\!\big(AU_V(-t)\varphi,BU_0(-t)\psi\big)\,dt\,
$$
can be established for arbitrary $\,\varphi\in{\cal H}_c\,$ and
$\,\psi\in{\cal H}.\,$ Further by the usage of intertwining
relationship $\,\,ZU_V(s)=U_0(s)Z\,\,$ and semigroup property
$\,\,U_V(t+s)=U_V(t)U_V(s)\,\,$ one gets
$$
\big(Z\varphi,\psi\big)\,\,-\,\,\big(U_0(s)U_V(-s)\varphi,\psi\big)\,\,=
\,\,-\,i\int_s^{\infty}\!\!\big(AU_V(-t)\varphi,BU_0(-t)\psi\big)\,dt\,
$$
and hence the inequality
$$
\big|\,\big(Z\varphi-U_0(s)U_V(-s)\varphi,\psi\big)\,\big|\,\,\leqslant\,\,C\,
\bigg(\int_s^{\infty}\!\!\|AU_V(-t)\varphi\|^2\,dt\bigg)^{1/2}\!\!\|\psi\|
$$
is valid with a certain constant $\,C>0\,$ which does not depend
on $\,\varphi\,$ and $\,\psi.\,$ Given $\,\varphi\in{\cal H}_c\,$
this implies the estimate
$$
\big\|\,Z\varphi-U_0(s)U_V(-s)\varphi\,\big\|\,\,\leqslant\,\,C\,
\bigg(\int_s^{\infty}\!\!\|AU_V(-t)\varphi\|^2\,dt\bigg)^{1/2},
$$
which extends to the whole $\,\cal H\,$ because $\,{\cal
H}_d\subset\ker Z\,$ and $\,U_V(t)=WU_0(t)Z.\,$ Finally, since
$\,AU_V(-t)\varphi\in{\rm L}_2(\mathbb R_+)\,$ for arbitrary
$\,\varphi\in{\cal H},\,$ we arrive at the conclusion that
$\,U_0(s)U_V(-s)\,$ converges strongly to $\,Z\,$ as
$\,s\to\infty.$

\bigskip
\medskip

\begin{center}
\bf \S\,6. One-dimensional Schroedinger operator
\end{center}

\bigskip

Proof of Theorem 3 reduces to verification of hypotheses imposed
in Theorem 1 in conformity with the Schroedinger operator
$\,\,L\,=\,T+V\,=\,-\,d^2/dx^2+V(x)\,$ acting in $\,{\cal H}={\rm
L}_2(\mathbb R_+)\,$ and determined in Section 2 by Dirichlet
boundary condition at zero. Henceforth bounded complex-valued
potential $\,V(x)\,$ is assumed to satisfy condition
(\ref{form4}).

\bigskip

{\bf Statement 5}\,\,(see\,[6]). {\it Resolvent
$\,R_V(\la)=(L-\lambda I)^{-1}$ of Schroedinger operator $\,L\,$
is meromorphic in $\,\mathbb C\setminus\mathbb R_+\,$ and its
poles coincide with eigenvalues of $\,L\,$ counting their
multiplicities. Discrete spectrum $\,\sigma_d(L)\,$ is finite
provided that operator $\,L\,$ has no spectral singularities,
while $\,\sigma_c(L)=\mathbb R_+.\,$ }

\bigskip

Note that for $\,\la=k^2,\,\,{\rm Im}\,k>0,\,$ the resolvent
$\,R_V(\la)\,$ of operator $\,L\,$ proves to be an integral
operator with the kernel
$$
R_V(x,\xi,\la)\,\,=\,\,\frac1{e(k)}\,\left\{
\begin{array}{rcl}s(x,k)\,e(\xi,k),&&x\leqslant\xi\\
e(x,k)\,s(\xi,k),&&\xi\leqslant x\\
\end{array}\right.
$$
where $\,s(x,k)\,$ is the solution of equation (\ref{form3})
satisfying initial conditions $\,s(0,k)=0,$ $\,s_x'(0,k)=1,\,$
while Jost solution $\,e(x,k)\,$ to equation (\ref{form3}) is
determined by asymptotics $\,e(x,k)\sim e^{ikx}\,$ at infinity.

\bigskip

{\bf Lemma 4} (cf. [16]). {\it For $\,\,k\in\mathbb C_+\,$ and
$\,\,x\in\mathbb R_+\,$ the following inequalities are valid
\begin{equation}
\big|s(x,k)\,e^{ikx}\big|\en\leqslant\en\min\left\{x,\frac1{|k|}\right\}\,\exp\left(\int_0^x\!\xi\,|V(\xi)|\,d\xi\right)\,,
\label{form6}
\end{equation}
\begin{equation}
\big|e(x,k)\,e^{-ikx}-1\big|\en\leqslant\en\exp\left(\int_x^{\infty}\!\!\min\left\{\xi,\frac1{|k|}\right\}\,|V(\xi)|\,d\xi\right)\,
-\,1\,. \label{form7}
\end{equation}
}

\medskip

In order to evaluate $\,s(x,k)\,$ note that it proves to be a
solution to integral equation
$$
s(x,k)\,\,\,=\,\,\,\frac{\sin kx}{k}\,\,\,+\,\,\int_0^x\frac{\sin
k(x-\xi)}{k}\,V(\xi)\,s(\xi,k)\,d\xi\,,
$$
whose integral kernel admits the estimate
$$
\left|\,\frac{\sin k(x-\xi)}{k}\,\right|\en\leqslant
\en\min\left\{x,\frac1{|k|}\right\}\,\exp\big({\rm
Im}\,k\,(x-\xi)\big)\,
$$
and consequently one has
$$
\big|s(x,k)\,e^{ikx}\big|\,\,\,\leqslant\,\,\,\min\left\{x,\frac1{|k|}\right\}
\left(1\,+\,\int_0^x|\,V(\xi)|\,|\,s(\xi,k)\,e^{i\xi
k}|\,d\xi\right).
$$
Estimate (\ref{form6}) now follows immediately by virtue of
Gronwall-Bellman inequality. As regards estimation of Jost
solution $\,e(x,k)\,$ one should take into account integral
equation
$$
e(x,k)\,\,\,=\,\,\,e^{ikx}\,\,\,-\,\,\int_x^{\infty}\frac{\sin
k(x-\xi)}{k}\,V(\xi)\,e(\xi,k)\,d\xi\,
$$
of Lippmann-Schwinger type which can be solved by iterations
method:
$$
e(x,k)\,e^{-ikx}\,=\,\,\sum_{n=0}^{\infty}\,\varepsilon^{(n)}(x,k)\,,\quad
\varepsilon^{(0)}(x,k)=1\,,
$$
$$\varepsilon^{(n)}(x,k)\,=\,\int_x^{\infty}\frac{e^{2ik(\xi-x)}-1}{2ik}\,V(\xi)\,
\varepsilon^{(n-1)}(\xi,k)\,d\xi\,.
$$
With regard for the inequality
$$
\,\,\displaystyle{\left|\,\frac{e^{2ik(\xi-x)}-1}{2ik}\,\right|\,\leqslant\,\min\left\{\xi,\frac1{|k|}\right\}},
$$
valid for $\,\xi\geqslant x,\,$ induction argument applies to
evaluate successive approximations
$$
|\,\varepsilon^{(n)}(x,k)|\,\,\,\leqslant\,\,\,\frac1{n!}\left(\int_x^{\infty}\min\left\{\xi,\frac1{|k|}\right\}
|V(\xi)|\,d\xi\right)^n\!
$$
and thus ensures estimate (\ref{form7}).

\bigskip

Denote by $\,A\,$ and $\,B\,$ operators of multiplication by
functions $\,a(x)\,$ and $\,b(x)\,$ to be chosen such that
\begin{eqnarray*}
\langle a
\rangle&:=&\left(\int_0^{\infty}x|\,a(x)|^2\,dx\right)^{1/2}<\,\infty\,,\\
\langle b
\rangle&:=&\left(\int_0^{\infty}x|\,b(x)|^2\,dx\right)^{1/2}<\,\infty\,.
\end{eqnarray*}

\medskip

{\bf Statement 6.} {\it Provided condition {\rm (\ref{form4})} is
satisfied operator function $\,e(k)AR_V(k^2)B^*\,$ extends
analytically to $\,\mathbb C_+\,$ and, moreover, for all
$\,k\in\mathbb C_+\,$ the inequality holds
$$
\big\|\,e(k)AR_V(k^2)B^*\,\big\|\,\,\,\leqslant\,\,\, K\,\langle a
\rangle\,\langle b \rangle\,,\quad K=\,\exp\,\langle
\sqrt{|\,V|}\,\rangle^2.
$$
}

\bigskip

{\bf Proof.} By virtue of (\ref{form6}) and (\ref{form7}) integral
kernel of the resolvent $\,R_V(\la),\,$ $\la=k^2,\,$ admits the
estimate
$$
\big|\,R_V(x,\xi,\la)\,\big|\,\,\,\leqslant\,\,\,\frac{K}{|\,e(k)|}\,\min\{x,\xi\}\,,\quad
K\,=\,\exp\left(\int_0^{\infty}\!x\,|V(x)|\,dx\right).
$$
As a consequence for arbitrary $\,f\in{\rm L}_2(\mathbb R_+)\,$
one has
\begin{multline*}
\|\,e(k)AR_V(k^2)B^*f\,\|^2\,\,\leqslant\,\,
K^2\int_0^{\infty}\!|\,a(x)|^2\,\bigg(\int_0^{\infty}\!\!\min\{x,\xi\}\,|\,b(\xi)|\,|f(\xi)|\,d\xi\,\bigg)^2dx\,\,\leqslant\\
\leqslant\,\,K^2\int_0^{\infty}\!\!\!|\,a(x)|^2\bigg(\int_0^{\infty}\!\!\!\big(\min\{x,\xi\}\big)^2|\,b(\xi)|^2\,d\xi\bigg)
\bigg(\int_0^{\infty}\!\!|f(\xi)|^2\,d\xi\bigg)\,dx\,\,\leqslant\\
\leqslant\,\,K^2\,\langle a \rangle^2\,\langle b \rangle^2\,
\|f\|^2\,.
\end{multline*}
Thus operator function $\,e(k)AR_V(k^2)B^*\,$ is analytic and
according to the above estimate bounded in the vicinity of each
pole of the resolvent $\,R_V(k^2)\,$ and therefore it has a
removable singularity therein.

\bigskip

{\bf Corollary} (cf. [1]). {\it For arbitrary $\,k\in\mathbb
C_+\,$ the inequality $\,\,\displaystyle{
\big\|AR_0(\lambda)A^*\big\|\,\leqslant\,\langle a \rangle^2}$ is
valid which implies smoothness of operator $\,A\,$ relative to
$\,T=T^*.$ }

\bigskip

Statements 5 and 6 guarantee that Theorem 1 applies to
Schroedinger operator
$$
\,\,L\,=\,T+V\,=\,-\,d^2/dx^2+V(x)\,\,
$$
with bounded potential $\,V(x)\,$ possessing finite first momentum
(\ref{form4}). In fact, the unperturbed operator
$\,T\,=\,-\,d^2/dx^2\,$ has purely continuous spectrum
$\,\sigma(T)=\mathbb R_+\,$ and according to Statement 5
conditions $\,(i)\,$ and $\,(ii)\,$ are satisfied provided that
operator $\,L\,$ has no spectral singularities.

Further, polar decomposition $\,V=J\,|V|,\,$ where
$\,|V|=\sqrt{V^*V},\,$ while $\,J={\rm sgn}\,V\,$ is a partial
isometry, produces an appropriate for our purposes factorization
$\,V=B^*A,\,$ so that $\,A=\sqrt{|V|}\,$ and
$\,B=\sqrt{|V|}\,J^*\,$ are operators of multiplication by
functions
$$
a(x)\,=\,\sqrt{|V(x)|}\,,\quad b(x)\,=\,\big(\,{\rm
sign}\,\overline{V(x)}\big)\,a(x)
$$
respectively, where $\,{\rm sign}\,z=z/|z|\,$ and $\,{\rm
sign}\,0=0.\,$ Under the condition (\ref{form4}) by virtue of
Statement 6 and the subsequent corollary operators $\,A\,$ and
$\,B\,$ are smooth relative to $\,T\,$ and, moreover, operator
function $\,Q_0(\la)=AR_0(\la)B^*\,$ is analytic and bounded in
$\,\mathbb C_\pm.\,$ If potential $\,V(x)\,$ decreases at infinity
at an appropriate rate $\,Q_0(\lambda)\,$ takes the values in the
trace class and Jost function $\,e(k)\,$ coincides (see. [17])
with Fredholm determinant
\begin{equation*}
e(k)\,\,=\,\,\det\big(I\,+\,AR_0(k^2)B^*\big).
\end{equation*}

Finally, due to the fact that Jost function $\,e(k)\,$ is analytic
in $\,\mathbb C_+,\,$ continuous up to the real axis and in virtue
of estimate (\ref{form7}) separated from zero at infinity,
Statement 6 implies condition $\,(iii)\,$ since and when operator
$\,L\,$ has no spectral singularities. Thus Schroedinger operator
$\,L=T+V\,$ under consideration satisfies all the hypotheses of
Theorem 1 and therefore its conclusion applies completing the
proof of Theorem 3.

\medskip

{\bf Corollary.} {\it Under the assumptions of Theorem {\rm 3}
condition $\,\sigma_d(L)=\varnothing\,$ gives a criterion of
similarity of operators $\,L\,$ and $\,T.\,$
}

\medskip

In conclusion we remark that Kato sufficient condition which
ensures similarity of operators $\,L\,$ and $\,T\,$ (see Section
1) is optimal in the sense that among the potentials $\,V(x)\,$
for which first momentum exceeds the critical value $1$ there
exist (see [10]) such that $\,\sigma_d(L)\ne\varnothing.\,$ The
estimates for the total number of eigenvalues and spectral
singularities of Schroedinger operator with complex potential were
obtained by the author in [13] and [18].

\bigskip
\medskip

\begin{center}
\bf References
\end{center}

\smallskip

\begin{enumerate}
\item Kato T. \, Wave operators and similarity for some
non-selfadjoint operators, Math. Annalen, 1966, v.162, p.258-279.

\item Schwartz J. \, Some non-selfadjoint operators, Comm. Pure
Appl. Math., 1960, v.13, p.609-639.

\item Jafaev D.R. \, Mathematical scattering theory, Transl. Math.
Monogr., v.105, Amer. Math. Soc., 1992.

\item Mochizuki K. \, On the large perturbation by a class of
non-selfadjoint operators, J. Math. Soc. Japan, 1967, v.19,
p.123-158.

\item Kako T., Yajima K. \, Spectral and scattering theory for a
class of non-selfadjoint operators, Sci. Papers of Coll. of Gen.
Educ. Univ. Tokyo, 1976, v.26, n.2, p.73-89.

\item Naimark M. A. \, Investigation of the spectrum and the
expansion in eigenfunctions of a non-selfadjoint differential
operator of the second order on a semi-axis, Proc. Mos. Math.
Soc., 1954, v.3, p.181-270.

\item Levin B. Ya. \, Transformations of Fourier and Laplace types
by means of solutions of differential equations of second order,
Doklady Math., 1956, v.106, n.2, p.187-190.

\item Glazman I. M.\, Direct methods of qualitative spectral
analysis of singular differential operators, Israel Prog.
Scientific Transl., 1965.

\item Hoffman K. \, Banach spaces of analytic functions,
Prentice-Hall, 1962.

\item Kato T. \, Perturbation theory for linear operators,
Springer-Verlag, 1966.

\item Marchenko V. A. \, Expansion in eigenfunctions of
non-selfadjoint singular second-order differential operators,
Sbornik Math., 1960, v.52, n.2, p.739-788.

\item Stepin S. A. \, Dissipative Schroedinger operator without a
singular continuous spectrum, Sbornik Math., 2004, v.195, n.6,
p.897-915.

\item Stepin S. A. \, On spectral components of the Schroedinger
operator with a complex potential, Russian Math. Surveys, 2013,
v.68, n.1, p.186-188.

\item Stankevich I. V. \, On linear similarity of certain
non-selfadjoint operators to selfadjoint operators and on the
asymptotic behavior for $t\to\infty$ of the solution of a
non-stationary Schroedinger equation, Sbornik Math., 1966, v.69,
n.2, p.161-207.

\item Keldysh M. V. \, On the completeness of the eigenfunctions
of some classes of non-selfadjoint linear operators, Russian Math.
Surveys, 1971, v.26, n.4, p.15-44.

\item Faddeev L. D. \, The inverse problem in the quantum theory
of scattering, Russian Math. Surveys, 1959, v.14, n.4, p.57-119.

\item Simon B. \, Resonances in one dimension and Fredholm
determinants, J. Funct. Anal., 2000, v.178, n.2, p.396-420.

\item Stepin S. A. \, Estimate for the number of eigenvalues of
the non-selfadjoint Schroedinger operator, Doklady Math., 2014,
v.89, n.2, p.202-205.
\end{enumerate}
\end{document}